\documentclass[11pt]{amsart}
\usepackage[a4paper, total={5.35 in, 8.3 in}]{geometry}
\usepackage{amsmath,amsfonts,amssymb}
\usepackage{amsfonts}
\usepackage{verbatim}
\usepackage{enumerate}
\usepackage{amsthm}
\usepackage{url}
\usepackage[colorlinks]{hyperref}
\usepackage[english]{babel}
\usepackage{tikz}
\newtheorem{theorem}{Theorem}[section]

\newtheorem{cor}[theorem]{Corollary}

\newcommand{\F}{\mathbb F}

\author[Lucas Reis]{Lucas Reis}
\address{Departamento de Matem\'{a}tica,
Universidade Federal de Minas Gerais,
UFMG,
Belo Horizonte MG (Brazil),
 30123-970}

\email{lucasreismat@gmail.com}

\title[A note on additive characters of finite fields]{A note on additive characters of finite fields}
\keywords{additive characters, finite fields, $\F_q$-order, reciprocal of polynomials}

\date{\today
}

\subjclass[2010]{12E20 (primary), 11T06(secondary)} 

\begin{document}
\maketitle

\begin{abstract}
Let $\F_q$ be the finite field with $q$ elements, where $q$ is a prime power and, for each integer $n\ge 1$, let $\F_{q^n}$ be the unique $n$-degree extension of $\F_q$. The $\F_q$-orders of an element in $\F_{q^n}$ and an additive character over $\F_{q^n}$ have been extensively used in the proof of existence results over finite fields (e.g., the Primitive Normal Basis Theorem). In this note we provide an interesting relation between these two objects.\end{abstract}

\section{Introduction}
Let $\F_q$ be the finite field with $q$ elements, where $q=p^s$ is a prime power and, for each integer $n\ge 1$, let $\F_{q^n}$ be the unique $n$-degree extension of $\F_q$. We observe that $\F_{q^n}$ can be viewed as an $\F_q$-vector space of dimension $n$. In this context, an element $\beta\in \F_{q^n}$ is normal over $\F_q$ if the set $\{\beta, \beta^q, \ldots, \beta^{q^{n-1}}\}$ of the conjugates of $\beta$ comprises an $\F_q$-basis for $\F_{q^n}$. The existence of normal elements is known for arbitrary $q$ and $n$ and, in fact, we have a closed formula for their number. If $\Phi_q$ denotes the Euler totient function for polynomials in $\F_q[x]$, there exist $\Phi_q(x^n-1)>0$ normal elements of $\F_{q^n}$ over $\F_q$. A proof of the latter can be obtained via the $\F_q$-order of elements in $\F_{q^n}$, which is defined as follows. The field $\F_{q^n}$ is endowed with the following $\F_{q}[x]$-module structure: for $g(x)=\sum_{i=0}^ma_ix^i\in \F_{q}[x]$ and $\alpha\in \F_{q^n}$, we set $g\circ \alpha=\sum_{i=0}^ma_i\alpha^{q^i}$. For an element $\alpha\in \F_{q^n}$, the set $\mathcal I_{\alpha}$ of the polynomials $g(x)\in \F_{q}[x]$ for which $g\circ \alpha=0$, is an ideal of $\F_q[x]$, hence principal.. The $\F_q$-order of $\alpha$ is the (unique) monic polynomial, denoted by $m_{\alpha, q}$, that generates the ideal $\mathcal I_{\alpha}$. Equivalently, $m_{\alpha, q}$ is the monic polynomial in $\F_{q}[x]$ of least degree such that $m_{\alpha, q}(x)\circ \alpha=0$. Since $(x^n-1)\circ \alpha=\alpha^{q^n}-\alpha$, we have that $x^n-1\in \mathcal I_{\alpha}$ for every $\alpha\in \F_{q^n}$; in particular, $m_{\alpha, q}$ divides $x^n-1$. The $\F_q$-order readily provides a criterion for normal elements: $\beta\in \F_{q^n}$ is normal over $\F_q$ if and only if $m_{\beta, q}(x)=x^n-1$. So normal elements are the ones of maximal $\F_q$-order. The formula $\Phi_q(x^n-1)$ for the number of normal elements is obtained through a simple inclusion-exclusion argument. In fact, a more general formula holds: if $f\in \F_q[x]$ is a monic divisor of $x^n-1$, the number of elements $\alpha\in \F_{q^n}$ with $m_{\alpha, q}=f$ equals $\Phi_q(f)$. For more details on these facts, see Theorem 11 of~\cite{ore} and Sec 2.2 of~\cite{thesis}.

The Primitive Normal Basis Theorem (PNBT) ensures the existence of normal elements $\beta\in \F_{q^n}$ that are primitive, i.e., generators of the cyclic (multiplicative) group $\F_{q^n}^*$. This result was first proved by Lenstra and Schoof~\cite{lenstra} and a proof without any use of computers was later given by Cohen and Huczynska~\cite{cohen}. The main idea is to use characters of $\F_{q^n}$ to build characteristic functions for the set of normal and primitive elements. This idea has been extensively used in the proof of the existence of elements in finite fields with many specified properties, beyond normality and primitivity. For more details see, Chapter 3 of~\cite{char} and the references therein. 

The characteristic function for the set of normal elements is obtained via additive characters. Write $q=p^s$, where $p$ is the characteristic of $\F_q$. An additive character of $\F_{q^n}$ is a function $\chi:\F_{q^n}\to \mathbb C^{\times}$ such that $\chi(a+b)=\chi(a)\cdot \chi(b)$. By the definition, each $a\in \F_{q^n}$ induces the additive character $\chi_a:\F_{q^n}\to \mathbb C^{\times}$ with
$$\chi_a(\alpha)=\exp\left({\frac{2\pi i \cdot \mathrm{Tr}_{q^n/p}(a\alpha)}{p}}\right),$$
where $\mathrm{Tr}_{q^n/p}(x)=\sum_{i=0}^{ns-1}x^{p^i}$ denotes the trace of $\F_{q^n}$ on $\F_p$. It is well known that every additive character of $\F_{q^n}$ is of such form. Moreover, from the identity of mappings $\chi_a\cdot \chi_b=\chi_{a+b}$, we have that the set $\widehat{\F_{q^n}}$ of additive characters of $\F_{q^n}$ is an abelian group (written multiplicatively), isomorphic to $\F_{q^n}$; the identity element is the trivial character $\chi_0$ with $\chi_0(\alpha)=1$ for every $\alpha\in \F_{q^n}$. 

The previous $\F_q[x]$-module structure of $\F_{q^n}$ lifts to the following $\F_q[x]$-module structure on $\widehat{\F_{q^n}}$: for $g\in \F_{q}[x]$ and $\chi \in \widehat{\F_{q^n}}$, we set $g\circ \chi:\alpha\to \chi(g\circ (\alpha))$, which is another element of $\widehat{\F_{q^n}}$. Within this structure, we have a natural extension of $\F_q$-order to additive characters: the set $\mathcal I_{\chi}$ of the polynomials $g(x)\in \F_{q}[x]$ for which $g\circ \chi=\chi_0$, is an ideal of $\F_q[x]$, hence principal. The $\F_q$-order of $\chi$ is the (unique) monic polynomial, denoted by $\mathrm{Ord}(\chi)$, that generates the ideal $\mathcal I_{\chi}$. It is direct to verify that, as before, $\mathrm{Ord}(\chi)\in \F_q[x]$ is a divisor of $x^n-1$.  In the proof of the PNBT~\cite{cohen, lenstra}, the character-sum formula for the characteristic function of normal elements depends on the sets $$\mathcal C_{f, q}:=\{\chi \in \widehat{\F_{q^n}}\,|\, \mathrm{Ord}(\chi)=f\},$$ where $f$ runs over the monic divisors of $x^n-1$ in $\F_q[x]$. We observe that the elements of $\mathcal C_{f, q}$ are of the form $\chi_a$ with $a\in \F_{q^n}$. Although no explicit description of the sets $C_{f, q}$ is required in the proof of the PNBT, further work needed to describe such sets for special values of $f$. For instance, see~\cite{panario}, where the cases $f(x)=1$ and $f(x)=x-1$ are considered, and the sets $\mathcal C_{1, q}=\{\chi_0\}$ and $\mathcal C_{x-1, q}=\{\chi_{a}\,|\, a\in \F_q^*\}$ are obtained.  The aim of this note is to provide the connection between the $\F_q$-orders of an element $\alpha\in \F_{q^n}$ and its associated additive character $\chi_{\alpha}\in \widehat{\F_{q^n}}$. Our result is stated as follows.
\begin{theorem}\label{thm:main}
Let $\alpha\in \F_{q^n}$ be an element of $\F_q$-order $f(x)=\sum_{i=0}^ma_ix^i$. Then the $\F_q$-order of $\chi_{\alpha}$ equals $f^{*}(x)=a_0^{-1}x^{m}f(1/x)$, the monic reciprocal of $f$.
\end{theorem}

As an immediate consequence of the previous theorem we have that, for each monic divisor $f\in \F_q[x]$ of $x^n-1$, the following holds:
$$\mathcal C_{f, q}=\{\chi_a\,|\, m_{a, q}=f^*\}.$$
In particular $\mathcal C_{1, q}=\{\chi_0\}$, $\mathcal C_{x-1, q}=\{\chi_{a}\,|\, a\in \F_q^*\}$ and $\mathcal C_{x^n-1, q}$ comprises the characters $\chi_{\beta}$ with $\beta\in \F_{q^n}$ a normal element over $\F_q$.

\section{Proof of Theorem~\ref{thm:main}}

We observe that, for monic polynomials $f, g\in \F_q[x]$ with $\gcd(f(x), x)=\gcd(g(x), x)=1$, we have the identity $(f^*)^*=f$ and $f$ divides $g$ if and only if $f^*$ divides $g^*$. Moreover, the $\F_q$-order of elements in $\F_{q^n}$ and $\widehat{\F_{q^n}}$ are divisors of $x^n-1$, hence relatively prime with $x$. From these facts it suffices to prove that, for every $\alpha\in \F_{q^n}$ and every monic $g\in \F_q[x]$ of degree at most $n-1$ with $\gcd(g(x), x)=1$, $g^*\circ \alpha=0$ if and only if $g\circ\chi_{\alpha}=\chi_0$. 

Write $q=p^s$ and $g(x)=a_mx^m+\sum_{i=0}^{m-1}a_ix^i\in \F_{q}[x]$, where $m<n$ and $a_m=1$. Notice that for every $a\in \F_{q^n}$, $$g\circ \chi_{\alpha}(a)=\exp\left(\frac{2\pi i\cdot \mathrm{Tr}_{q^n/p}(\alpha\cdot L_g(a))}{p}\right),$$ where $L_g(x)=x^{q^m}+\sum_{i=0}^{m-1}a_ix^{q^i}$. Therefore, it suffices to prove that $g^{*}\circ \alpha=0$ if and only if $$\mathrm{Tr}_{q^n/p}(\alpha L_g(x))\equiv 0\pmod {x^{q^n}-x}.$$ Let $M(x)\in \F_q[x]$ be the unique polynomial of degree at most $q^n-1$ such that $\mathrm{Tr}_{q^n/p}(\alpha L_g(x))\equiv M(x)\pmod {x^{q^n}-x}$. It is clear that $M(x)$ is of the form $\sum_{i=0}^{ns-1}c_ix^{p^i}$. Since $\mathrm{Tr}_{q^n/p}(a)\in \F_p$ for every $a\in \F_{q^n}$, we have that $M(x)^p\equiv M(x)\pmod {x^{q^n}-x}$. In particular, $c_j=c_{j-1}^p$ for every $1\le j\le n$. Therefore, $M(x)$ vanishes if and only if $c_0=0$. Let us compute the coefficient $c_0$. We observe that 
$$\mathrm{Tr}_{q^n/p}(\alpha L_g(x))=\sum_{i=0}^{ns-1}\left(\alpha\sum_{j=0}^{m-1}a_jx^{q^j}\right)^{p^i}=\sum_{j=0}^{m-1}\sum_{i=0}^{ns-1}a_j^{p^i}\alpha^{p^i}x^{p^{sj+i}}.$$

Since $x^{p^{sn}}\equiv x\pmod {x^{q^n}-x}$, in the last sum, the terms contributing to $c_0$ are the ones of the form $a_j^{p^i}\alpha^{p^i}x^{p^{i+sj}}$ with $i+sj\equiv 0\pmod {sn}$, where $0\le i\le sn-1$ and $0\le j\le m<n$. The latter implies that $i+sj\equiv 0\pmod {sn}$ if and only if $i=j=0$ or $i+sj=sn$, i.e., $(i, j)=(s(n-t), t)$ with $0\le t\le m$.  The corresponding terms are $a_{t}^{p^{s(n-t)}}\alpha^{p^{s(n-t)}}=a_{t}\alpha^{q^{n-t}}$, since $a_t\in \F_{q}=\F_{p^s}$ for every $1\le t\le m$. Therefore, we have the following equality:
$$c_0=a_0\alpha+a_1\alpha^{q^{n-1}}+\ldots+a_m\alpha^{q^{n-m}}.$$
Raising both sides of the previous equality to the $q^m$-th power, we obtain that $$c_0^{q^m}=a_0\alpha^{q^m}+\ldots+a_m\alpha=a_0\cdot (g^{*}\circ \alpha).$$
Since $a_0=g(0)\ne 0$, we have that $c_0=0$ if and only if $g^{*}\circ \alpha=0$, concluding the proof.

\section{Additional Remarks}

A $t$-degree monic polynomial $f\in \F_q[x]$ is self-reciprocal if it coincides with its monic reciprocal, i.e., $$f^*(x)=f(0)^{-1}x^tf(1/x)=f(x).$$ The following corollary is an immediate consequence of Theorem~\ref{thm:main}.

\begin{cor}\label{cor:srim}
The $\F_q$-orders of $\alpha\in \F_{q^n}$ and its associated additive character $\chi_{\alpha}\in \widehat{\F_{q^n}}$ coincide if and only if one of them is self-reciprocal.
\end{cor}

In Subsec.~2.3 of~\cite{KR} and Subsec.~2.3.1 of~\cite{R}, the authors mistakenly claimed that the $\F_q$-orders of $\alpha$ and its associated character $\chi_{\alpha}$ coincide for every $\alpha\in \F_{q^n}$. Fortunately, the results on these papers are not affected by this false claim. As in the proof of the PNBT, the results of~\cite{R} do not require explicit description of sets $\mathcal C_{f, q}$ (see Subsec. 3.1 of~\cite{R}). In particular, the claim relating the $\F_q$-orders of $\alpha$ and $\chi_{\alpha}$ are not employed in~\cite{R}. Moreover, in~\cite{KR}, this claim is implicitly used only for characters of $\F_q$-order $1$ and $(x-1)^m$, where $m$ is a power of the characteristic $p$ (see Theorem 3.1~in~\cite{KR}). The latter comprise only self-reciprocal polynomials, where the claim is true by Corollary~\ref{cor:srim}.

Applying the previous corollary, we characterize the positive integers $n$ for which the false claim in~\cite{KR, R} holds true. We have seen that the $\F_q$-orders of $\alpha$ and $\chi_{\alpha}$ are divisors of $x^n-1$. Conversely, for each monic divisor $f\in \F_q[x]$ of $x^n-1$, there exist elements in $\F_{q^n}$ with $\F_q$-order $f$. Hence the $\F_q$-orders of $\alpha$ and $\chi_{\alpha}$ coincide for every $\alpha\in \F_{q^n}$ if and only if every monic divisor of $x^n-1$ over $\F_q$ is a self-reciprocal polynomial. Equivalently, every monic irreducible divisor of $x^v-1$ over $\F_q$ is self-reciprocal, where $n=p^u\cdot v$ with $\gcd(v, p)=1$. According to Theorem~1 of~\cite{meyn}, the latter holds if and only if $q^j\equiv -1\pmod v$ for some positive integer $j\le v$. We obtain the following corollary.

\begin{cor}
Write $n=p^u\cdot v$, where $\gcd(v, p)=1$. The $\F_q$-orders of $\alpha$ and $\chi_{\alpha}$ coincide for every $\alpha\in \F_{q^n}$ if and only if $q^j\equiv -1\pmod v$ for some $1\le j\le v$.
\end{cor}


\end{document}